\documentclass[12pt] {article}
\usepackage{amsmath,amsthm}
\usepackage{amssymb,latexsym}
\usepackage{amscd}
\title{Solvability of the quaternionic Monge-Amp\`ere
equation on compact manifolds with a flat hyperK\"ahler metric.}
\date{}
\author{Semyon Alesker \footnote{Partially supported by ISF grant 701/08.}
\\  { \normalsize Department of Mathematics, Tel Aviv University, Ramat Aviv}
\\  { \normalsize 69978 Tel Aviv, Israel }
\\ {\normalsize e-mail: semyon@post.tau.ac.il}}

\def\eps{\varepsilon}
\def\alp{\alpha}
\def\ome{\omega}
\def\Ome{\Omega}
\def\lam{\lambda}
\def\Lam{\Lambda}

\def\to{\rightarrow}
\def\qed { Q.E.D. }
\def\pt{\partial}

\def\RR{\mathbb{R}}
\def\CC{\mathbb{C}}

\def\HH{\mathbb{H}}

\swapnumbers
\newtheorem{theorem}{Theorem}[section]
\newtheorem{corollary}[theorem]{Corollary}
\newtheorem{lemma}[theorem]{Lemma}
\newtheorem{proposition}[theorem]{Proposition}
\newtheorem{claim}[theorem]{Claim}
\theoremstyle{definition}
\newtheorem{example}[theorem]{Example}
\newtheorem{definition}[theorem]{Definition}
\newtheorem{remark}[theorem]{Remark}

\theoremstyle{proposition-definition}
\newtheorem{proposition-definition}[theorem]{Proposition-Definition}

\numberwithin{equation}{section}

 \def\cb{{\cal B}} \def\cc{{\cal C}}
\def\cd{{\cal D}}  
 \def\ch{{\cal H}} 
 \def\ck{{\cal K}} 
\def\cm{{\cal M}} \def\cn{{\cal N}} \def\co{{\cal O}}

  \def\cu{{\cal U}}
\def\cv{{\cal V}}

\begin{document}
\maketitle

\begin{abstract}
A quaternionic version of the Calabi problem  was formulated in
\cite{alesker-verbitsky-10}. It conjectures a solvability of a
quaternionic Monge-Amp\`ere equation on a compact HKT manifold (HKT
stays for HyperK\"ahler with Torsion). In this paper this problem is
solved under the extra assumption that the manifold admits a flat
hyperK\"ahler metric compatible with the underlying hypercomplex
structure. The proof uses the continuity method and {\itshape a
priori} estimates.
\end{abstract}

\tableofcontents \setcounter{section}{-1}

\section{Introduction.}\label{S:introduction}
In recent years there was suggested a quaternionic analogue of the
classical real and complex Monge-Amp\`ere equations. Thus in
\cite{alesker-jga-03} the author has introduced quaternionic
Monge-Amp\`ere equation on the flat quaternionic space $\HH^n$ and
solved the Dirichlet problem for it under appropriate assumptions.
Then M. Verbitsky and the author \cite{alesker-verbitsky-06} have
generalized the equation to the broader class of so called
hypercomplex manifolds. They have also formulated a conjecture about
existence of a solution of this quaternionic Monge-Amp\`ere equation
which is a quaternionic analogue of the well known Calabi problem
for the complex Monge-Amp\`ere equation. Moreover they proved a
uniform a priori estimate for its solution (under some assumption)
and uniqueness of a solution up to a constant. The goal of this
paper is to solve the conjecture in the case of a compact
hypercomplex manifold which admits a flat hyperK\"ahler metric
compatible with the underlying hypercomplex structure (see Theorem
\ref{T:main-theorem} below).

Recall that the original Calabi problem for the complex
Monge-Amp\`ere equation was formulated by him in 1954. It was
eventually solved by Yau in 1976 \cite{yau}. Before this, Aubin
\cite{aubin-calabi} had made significant progress towards its proof.
A real version of the Calabi problem was formulated and solved by
Cheng and Yau \cite{cheng-yau}.

Let us also mention that recently Harvey and Lawson
\cite{harvey-lawson-09} have extended the notion of (homogeneous)
Monge-Amp\`ere equation beyond real, complex, and quaternionic
cases.

\hfill

In order to formulate the main result precisely, let us recall the
notions of hypercomplex and HKT-manifolds.

\begin{definition}
A {\itshape hypercomplex} manifold is a smooth manifold $M$ together
with a triple $(I,J,K)$ of complex structures satisfying the usual
quaternionic relations:
$$IJ=-JI=K.$$
\end{definition}
\begin{remark}
(1) We suppose here that the complex structures $I,J,K$ act on the
{\itshape right} on the tangent bundle $TM$ of $M$. This action
extends uniquely to the right action of the algebra $\HH$ of
quaternions on $TM$.

(2) It follows that the dimension of a hypercomplex manifold $M$ is
divisible by 4.

(3) Hypercomplex manifolds were explicitly introduced by Boyer
\cite{boyer}.
\end{remark}

Let $(M^{4n},I,J,K)$ be a hypercomplex manifold. Let us denote by
$\Lam^{p,q}_I(M)$ the vector bundle of differential forms of the
type $(p,q)$ on the complex manifold $(M,I)$. By the abuse of
notation we will also denote by the same symbol $\Lam^{p,q}_{I}(M)$
the space of $C^\infty$-sections of this bundle.

Let
\begin{eqnarray}\label{l1}
\pt\colon \Lambda_I^{p,q}(M)\to \Lambda_I^{p+1,q}(M)
\end{eqnarray}
be the usual $\pt$-differential on differential forms on the complex
manifold $(M,I)$.

Set
\begin{eqnarray}\label{l2}
\pt_J:=J^{-1}\circ \bar \pt \circ J.
\end{eqnarray}

\begin{claim}[\cite{verbitsky-hkt}]\label{l3}
(1)$ J\colon\Lambda_I^{p,q}(M)\to\Lam_I^{q,p}(M).$

(2) $\pt_J\colon \Lambda_I^{p,q}(M)\to\Lam_I^{p+1,q}(M).$

(3) $\pt\pt_J=-\pt_J\pt$.
\end{claim}

\begin{definition}[\cite{verbitsky-hkt}]\label{l4}
Let $k=0,1,\dots,n$. A form $\ome\in \Lam^{2k,0}_I(M)$ is called
\itshape{real} if
$$\overline{J\circ \ome}=\ome.$$
\end{definition}

We will denote the subspace of real $C^\infty$-smooth $(2k,0)$-forms
on $(M,I)$ by $\Lam^{2k,0}_{I,\RR}(M)$.
\begin{lemma}\label{l5}
Let $(M,I,J,K)$ be a hypercomplex manifold. Let $f\colon M\to \RR$
be a smooth function. Then $\pt\pt_J f\in \Lam^{2,0}_{I,\RR}(M)$.
\end{lemma}

We call $\pt\pt_J h$ the quaternionic Hessian of $f$. In many
respects it is analogous to the usual real and complex Hessians. It
becomes particularly transparent on the flat space $\HH^n$ where it
can be written in coordinates; see the discussion in Section
\ref{S:HKT-manifolds}.

\begin{definition}\label{Def:positive-forms}
Let $\ome\in \Lam^{2,0}_{I,\RR}(M)$. Let us say that $\ome$ is
non-negative (notation: $\ome\geq 0$) if
$$\ome(Y,Y\circ J)\geq 0$$
for any (real) vector field $Y$ on the manifold $M$. The form $\ome$
is called strictly positive (notation: $\ome>0$) if $\ome(Y,Y\circ
J)>0$ for any non-vanishing (real) vector field $Y$.

Equivalently, $\omega$ is non-negative (resp., strictly positive) if
and only if $\omega(Z, \bar Z \circ J)\geq 0$ (resp., $>0$) for any
non-vanishing $(1,0)$-vector field $Z$.
\end{definition}

Let $g$ be a Riemannian metric on a hypercomplex manifold $M$. The
metric $g$ is called {\itshape quaternionic Hermitian} (or
hyperhermitian) if $g$ is invariant with respect to the group
$SU(2)\subset \HH$ of unitary quaternions, i.e. $g(X\cdot q,Y\cdot
q)=g(X,Y)$ for any (real) vector fields $X,Y$ and any $q\in \HH$
with $|q|=1$.

Given a quaternionic Hermitian metric $g$ on a hypercomplex manifold
$M$, consider the differential form
$$\Ome:=\ome_J-\sqrt{-1}\ome_K$$
where $\ome_L(A,B):=g(A,B\circ L)$ for any $L\in \HH$ with $L^2=-1$,
and any real vector fields $A,B$ on $M$. It is easy to see that
$\Ome$ is a $(2,0)$-form with respect to the complex structure $I$.
Moreover $\Ome$ is real in the sense of Definition
\ref{Def:positive-forms}, thus $\Ome\in \Lam^{2,0}_{I,\RR}(M)$.

\begin{definition}\label{def-hkt-metr}
The metric $g$ on $M$ is called HKT-metric if
$$\pt \Ome =0.$$
\end{definition}
We call such a form $\Ome$, corresponding to an HKT-metric, an
HKT-form.

\begin{remark}
HKT manifolds were introduced in the physical literature by Howe and
Papadopoulos \cite{howe-papa}. For the mathematical treatment see
Grantcharov-Poon \cite{grantcharov-poon} and Verbitsky
\cite{verbitsky-hkt}. The original definition of HKT-metrics in
\cite{howe-papa} was different but equivalent to Definition
\ref{def-hkt-metr}; the latter was given in \cite{grantcharov-poon}.
\end{remark}

\begin{remark}\label{R:hyperKahler}
The classical hyperK\"ahler metrics (i.e. Riemannian metrics with
holonomy contained in the group $Sp(n)$) form a subclass of
HKT-metrics. It is well known that a quaternionic Hermitian metric
$g$ is hyperK\"ahler if and only if the form $\Ome$ is closed, or
equivalently $\pt\Ome=\bar\pt\Ome=0$.
\end{remark}

Now we can formulate the main result.
\begin{theorem}\label{T:main-theorem}
Let $(M^{4n},I,J,K)$ be a compact connected hypercomplex manifold
with an HKT form $\Ome_0$.\footnote{In this paper all HKT-metrics,
and consequently HKT-forms, are assumed to be infinitely smooth.}
Let us assume in addition that it admits a flat hyperK\"ahler metric
compatible with the underlying hypercomplex structure. Let $f\in
C^\infty(M)$ be a real valued function. Then there exists a unique
constant $A$ such that the quaternionic Monge-Amp\`ere equation
\begin{eqnarray}\label{E:main-equation}
(\Ome_0+\pt\pt_J\phi)^n=A e^f\Ome_0^n
\end{eqnarray}
has a $C^\infty$-smooth solution.
\end{theorem}
\begin{remark}
(1) It was shown in \cite{alesker-verbitsky-10} that solution $\phi$
is unique up to an additive constant.

(2) The constant $A$ is determined as follows. Let $\Ome$ be the
HKT-form corresponding to the flat hyperK\"ahler metric whose
existence is assumed in the theorem. Then $A$ is found from the
equation
$$\int_M A e^f\cdot \Ome_0^n\wedge\bar\Ome^n=\int_M
\Ome_0^n\wedge\bar\Ome^n.$$

(3) This theorem was conjectured by M. Verbitsky and the author in
\cite{alesker-verbitsky-10} in a more general form: without the
assumption of existence a flat hyperK\"ahler metric.

(4) The equation (\ref{E:main-equation}) is a non-linear second
order elliptic differential equation. The ellipticity was shown in
\cite{alesker-verbitsky-10}.

(5) Existence of a flat hyperK\"ahler metric implies that the
hypercomplex structure $(I,J,K)$ is locally flat, i.e. locally
isomorphic to the standard flat space $\HH^n$ of $n$-tuple of
quaternions.

(6) This theorem can be stated in a slightly more refined form
involving H\"older spaces rather than $C^\infty$, see Theorem
\ref{T:MA-solution} below.

(7) The obvious example of a hypercomplex manifold $M$ satisfying
the assumptions of the theorem is a quaternionic torus: quotient of
$\HH^n$ by a lattice. However there are more examples coming from
the Bieberbach classification of crystallographic groups (see e.g.
\cite{vinberg-shvartsman}).
\end{remark}

Notice that recently Verbitsky \cite{verbitsky-balanced-09} has
suggested a geometric interpretation of solutions of the equation
(\ref{E:main-equation}) under appropriate assumptions on the right
hand side.

The proof of the theorem uses the continuity method and a priori
estimates. The standard elliptic regularity machinery, discussed in
Section \ref{S:higher-order}, implies that it suffices to prove a
$C^{2,\alpha}$ a priori estimate for some $\alp\in (0,1)$.  The
$C^0$ estimate was obtained first in \cite{alesker-verbitsky-10}
under more general assumptions than in Theorem \ref{T:main-theorem}.
Very recently Shelukhin and the author \cite{alesker-shelukhin} have
obtained a $C^0$ estimate by a different method and under different
assumptions than in \cite{alesker-verbitsky-10} which however are
also satisfied in Theorem \ref{T:main-theorem}. The main point of
this paper is to make two following steps: first to obtain $C^0$
estimate on a Laplacian of $\phi$ (Section \ref{S:second_order}),
and then to deduce from it a $C^{2,\alp}$ estimate (Section
\ref{S:2-alpha}). The first step uses a modification of the well
known Pogorelov's method. This modification is not completely
straightforward, and this is exactly the step where all the
assumptions of the theorem are used, i.e. existence of a flat
hyperK\"ahler metric. The second step uses a quaternionic version of
the Evans-Krylov method (see Section \ref{S:2-alpha} for further
references). It works under more general assumptions, namely on
manifolds with locally flat hypercomplex structure (which may not
admit a compatible hyperK\"ahler metric).

In Section \ref{S:linear-algebra} we recall relevant definitions and
facts from the quaternionic linear algebra. In Section
\ref{S:HKT-manifolds} we recall few facts on HKT-manifolds. These
two sections contain no new results, they are added for convenience
of the reader only.

\hfill

{\bf Acknowledgement.} I thank M. Verbitsky for numerous very useful
discussions.

\section{Quaternionic linear algebra.}\label{S:linear-algebra}
The standard theory of vector spaces, basis, and dimension works
over any non-commutative field, e.g. $\HH$, exactly like in the
commutative case. The only remark is that one should distinguish
between right and left vector spaces. The two cases are completely
parallel. We will restrict to the case of right vector spaces, i.e.
vectors are multiplied by scalars on the right.

However the theory of non-commutative determinants is quite
different and deserves special discussion. We will need to remind
the notion of Moore determinant on the class of quaternionic
matrices called hyperhermitian.\footnote{The Moore determinant was
used in the original approach of \cite{alesker-bsm} to define the
quaternionic Monge-Amp\`ere operator on the flat space $\HH^n$.
Later on, this operator was generalized in
\cite{alesker-verbitsky-06} to more general class of hypercomplex
manifolds without using explicitly the Moore determinant. However
this notion often still seems to be convenient while working on the
flat space; in particular it will be used extensively in this
paper.} They are analogues of real symmetric and complex hermitian
matrices. The behavior of the Moore determinant of such matrices is
analogous in many respects to the behavior of the usual determinant
of real symmetric and complex hermitian matrices. We believe that
any general identity or inequality known for usual determinant of
the real symmetric or complex hermitian matrices can be generalized
to the Moore determinant of hyperhermitian matrices, though the
proofs might be slightly more tricky. Here we review some of the
relevant material. The discussion on determinants follows
\cite{alesker-bsm} where most of the proofs can be found. Another
good reference to quaternionic determinants is \cite{aslaksen}; for
a relation of quaternionic determinants to a general theory
\cite{gelfand-retakh-wilson-survey} of non-commutative (quasi-)
determinants see \cite{gelfand-retakh-wilson}.
\begin{definition}\label{det-5}
Let $V$ be a right $\HH$-vector space.  A {\itshape hyperhermitian
semilinear form} on $V$ is a map $ a:V \times V \to \HH$ satisfying
the following properties:

(a) $a$ is additive with respect to each argument;

(b) $a(x,y \cdot q)= a(x, y) \cdot q$ for any $x,y \in V$ and any
$q\in \HH$;

(c) $a(x,y)= \overline{a(y,x)}$.
\end{definition}

\begin{remark}\label{R:hyper-forms}
Hyperhermitian semi-linear forms on $V$ are in bijective
correspondence with real valued quadratic forms on the underlying
real space ${}\!^ {\RR} V $ of $V$
$$b\colon {}\!^ {\RR} V\to \RR$$
which are invariant under multiplication by the norm one
quaternions, i.e. $b(x\cdot q)=b(x)$ for any $x\in V$ and any $q\in
\HH$ with $|q|=1$.
\end{remark}

\begin{example}\label{det-6} Let $V= \HH ^n$ be the standard coordinate space
considered as right vector space over $\HH$. Fix a {\itshape
hyperhermitian} $n \times n$-matrix $(a_{ij})_{i,j=1}^{n}$, i.e.
$a_{ij} =\bar a_{ji}$, where $\bar q$ denotes the usual quaternionic
conjugation of $q\in \HH$. For $x=(x_1, \dots, x_n), \, y=(y_1,
\dots, y_n)$ define
$$A(x,y) = \sum _{i,j} \bar x_i a_{ij} y_j$$
(notice the order of the terms!). Then $A$ defines hyperhermitian
semilinear form on $V$.

The set of all hyperhermitian $n\times n$-matrices will be denoted
by $\ch_n$. Then $\ch_n$ a vector space over $\RR$.
\end{example}
In general one has the following standard claims.

\begin{claim}\label{det-7} Fix a basis in a finite dimensional right quaternionic
vector space $V$. Then there is a natural bijection between the
space of hyperhermitian semilinear forms on $V$ and the space
$\ch_n$ of $n \times n$-hyperhermitian matrices.
\end{claim}
This bijection is in fact described in previous Example \ref{det-6}.

\begin{claim}\label{det-8} Let $A$ be a matrix
of the given hyperhermitian form in a given basis. Let $C$ be
transition matrix from this basis to another one. Then the matrix
$A'$ of the given form in the new basis is equal $$A' =C^* AC,$$
where $(C^*)_{ij}=\bar C_{ji}$.
\end{claim}
\begin{remark}\label{det-9}
Note that for any hyperhermitian matrix $A$ and for any matrix $C$
the matrix $C^* AC$ is also hyperhermitian. In particular the matrix
$C^* C$ is always hyperhermitian.
\end{remark}

\begin{definition}\label{det-10} A hyperhermitian semilinear form $a$
is called {\itshape positive definite} if $a(x,x)>0$ for any
non-zero vector $x$. Similarly $a$ is called {\itshape non-negative
definite} if $a(x,x)\geq 0$ for any vector $x$.
 \end{definition}

Let us fix on our quaternionic right vector space $V$ a positive
definite hyperhermitian form $( \cdot , \cdot )$. The space with
fixed such a form will be called {\itshape hyperhermitian} space.

For any quaternionic linear operator $\phi: V\to V$ in
hyperhermitian space one can define the adjoint operator $\phi ^* :V
\to V$ in the usual way, i.e. $(\phi x,y)= (x, \phi ^* y)$ for any
$x,y \in V$. Then if one fixes an orthonormal basis in the space $V$
then the operator $\phi$ is selfadjoint if and only if its matrix in
this basis is hyperhermitian.

\begin{claim}\label{det-11}
For any  selfadjoint operator in a hyperhermitian space there exists
an orthonormal basis such that its matrix in this basis is diagonal
and real.
\end{claim}
 Now we are going to define the Moore
determinant of  hyperhermitian matrices. The definition below is
different from the original one \cite{moore} but equivalent to it.

Any quaternionic matrix $A\in M_n(\HH)$ can be considered as a
matrix of an $\HH$-linear endomorphism of $\HH^n$. Identifying
$\HH^n$ with $\RR^{4n}$ in the standard way we get an $\RR$-linear
endomorphism of $\RR^{4n}$. Its matrix in the standard basis will be
denoted by ${}^{\mathbb{R}} A$, and it is called the realization of
$A$. Thus ${}^{\mathbb{R}} A\in M_{4n}(\RR)$.

 Let us consider the entries of $A$ as formal variables (each
quaternionic entry corresponds to four commuting real variables).
Then $\det ({}^{\mathbb{R}} A)$  is a homogeneous polynomial of
degree $4n$ in $n(2n-1)$ real variables. Let us denote by $Id$ the
identity matrix.
 One has the following result.
\begin{theorem}\label{det-12}
There exists a polynomial $P$ defined on the space $\ch_n$ of all
hyperhermitian $n \times n$-matrices such that for any
hyperhermitian $n \times n$-matrix $A$ one has $\det({}^{\mathbb{R}}
A)= P^4(A)$ and $P(Id)=1$. $P$ is defined uniquely by these two
properties. Furthermore $P$ is homogeneous of degree $n$ and has
integer coefficients.
\end{theorem}
Thus for any hyperhermitian matrix $A$ the value $P(A)$ is a real
number, and it is called the {\itshape Moore determinant} of the
matrix $A$. The explicit formula for the Moore determinant  was
given by Moore \cite{moore} (see also \cite{aslaksen}). From now on
the Moore determinant of a matrix $A$ will be denoted by $\det A$.
This notation should not cause any confusion with the usual
determinant of real or complex matrices due to part (i) of the next
theorem.
\begin{theorem}\label{det-13}

(i) The Moore determinant of any complex hermitian matrix considered
as quaternionic hyperhermitian matrix is equal to its usual
determinant.

(ii) For any hyperhermitian $n\times n$-matrix $A$ and any matrix
$C\in M_n(\HH)$ the Moore determinant satisfies
$$\det (C^*AC)= \det A \cdot \det(C^*C).$$
\end{theorem}
\begin{example}\label{det-14}

(a) Let $A =diag(\lam_1, \dots, \lam _n)$ be a diagonal matrix with
real $\lam _i$'s. Then $A$ is hyperhermitian and the Moore
determinant $\det A= \prod _{i=1}^n \lam_i$.

(b)  A general hyperhermitian $2 \times 2$-matrix $A$ has the form
 $$ A=  \left[ \begin {array}{cc}
                     a&q\\
                \bar q&b\\
                \end{array} \right] ,$$
where $a,b \in \RR, \, q \in \HH$. Then $\det A =ab - q \bar q$.
\end{example}

\begin{definition}\label{det-14.5}
A hyperhermitian $n\times n$-matrix $A=(a_{ij})$ is called {\itshape
positive} (resp. {\itshape non-negative}) {\itshape definite} if for
any non-zero vector $\xi= \left[\begin{array}{c}
             \xi_1\\
             \vdots\\
             \xi_n
             \end{array}\right]$ one has
$\xi^*A\xi=\sum_{ij}\bar\xi_ia_{ij}\xi_j
>0$ (resp. $\geq 0$).
\end{definition}

\begin{claim}\label{det-15}
Let $A$ be  a non-negative (resp. positive) definite hyperhermitian
matrix. Then $\det A \geq 0 \, (\mbox{ resp. } \det A
>0)$.
\end{claim}

Moreover there is a version of the Sylvester criterion of positive
definiteness of a hyperhermitian matrix. It is formulated in terms
of the Moore determinants and is completely analogous to the
classical real and complex results, see \cite{alesker-bsm}, Theorem
1.1.13.

Let us remind now the definition of the mixed determinant of
hyperhermitian matrices in analogy with the case of real symmetric
matrices \cite{aleksandrov-38}.
\begin{definition}\label{det-16}
Let $A_1, \dots ,A_n$ be hyperhermitian $n \times n$- matrices.
Consider the homogeneous polynomial in real variables $\lam _1
,\dots , \lam _n$ of degree $n$ equal to $\det(\lam_1 A_1 + \dots +
\lam_n A_n)$. The coefficient of the monomial $\lam_1 \cdot \dots
\cdot \lam_n$ divided by $n!$ is called the {\itshape mixed
determinant} of the matrices  $A_1, \dots ,A_n$, and it is denoted
by $\det(A_1, \dots ,A_n)$.
\end{definition}
Note that the mixed determinant is symmetric with respect to all
variables, and linear with respect to each of them, i.e.
$$\det (\lam A_1' +\mu A_1'', A_2, \dots, A_n )=
\lam \cdot \det( A_1', A_2, \dots, A_n ) + \mu \cdot \det(A_1'',
A_2, \dots, A_n )$$ for any {\itshape real} $\lam , \, \mu$. Note
also that $\det(A, \dots, A)=\det A$.

\begin{theorem}\label{det-17}
The mixed determinant of positive (resp. non-negative) definite
matrices is positive (resp. non-negative).
\end{theorem}
This theorem is proved in \cite{alesker-bsm}, Theorem 1.1.15(i).
Moreover a version of the A.D. Aleksandrov inequality for mixed
determinants can be proven, see \cite{alesker-bsm}, Theorem 1.1.15
and Corollary 1.1.16.

\section{HKT manifolds.}\label{S:HKT-manifolds}
In this section we recall few facts about HKT-manifolds in addition
to those stated in the introduction.

\begin{definition}[\cite{alesker-verbitsky-06}]\label{l6}
Let $(M,I,J,K)$ be a hypercomplex manifold. A $C^2$-smooth function
$$h:M\to \RR$$ is called quaternionic plurisubharmonic if $\pt\pt_J
h$ is a non-negative section of $\Lambda^{2,0}_{I,\RR}(M)$. $h$ is
called strictly plurisubharmonic if $\pt\pt_J h$ is strictly
positive at every point.
\end{definition}

\begin{remark}
The notion of quaternionic plurisubharmonicity can be generalized to
continuous functions, see \cite{alesker-verbitsky-06}, Section 5. On
the flat space $\HH^n$ this notion was earlier defined even for
upper semi-continuous functions in \cite{alesker-bsm}.
\end{remark}

Let us discuss the relations of plurisubharmonic functions to the
HKT-geometry. Let us denote by $S_\HH(M)$ the vector bundle over a
hypercomplex manifold $M$ such that its fiber over a point $x\in M$
is equal to the space of hyperhermitian forms on the tangent space
$T_xM$. Consider the map of vector bundles
\begin{eqnarray}\label{E:isomor-t}
t\colon \Lam^{2,0}_{I,\RR}(M)\to S_\HH(M)
\end{eqnarray}
defined by $t(\eta)(A,A)=\eta(A,A\circ J)$ for any (real) vector
field $A$ on $M$. Then $t$ is an isomorphism of vector bundles (this
was proved in \cite{verbitsky-hkt}).

\begin{theorem}[\cite{alesker-verbitsky-06}, Prop. 1.14]\label{m}
(1) Let $f$ be an infinitely smooth strictly plurisubharmonic
function on a hypercomplex manifold $(M,I,J,K)$. Then $t(\pt\pt_J
f)$ is an HKT-metric.

(2) Conversely assume that $g$ is an HKT-metric. Then any point
$x\in M$ has a neighborhood $U$ and an infinitely smooth strictly
plurisubharmonic function $f$ on $U$ such that $g=t(\pt\pt_J f)$ in
$U$. Equivalently $\Ome=\pt\pt_J f$, where $\Ome$ is the HKT-form
corresponding to $g$ (as defined in the introduction).
\end{theorem}

In this paper we will often work with the flat hypercomplex manifold
$\HH^n$. In this case there is an equivalent way to rewrite the
quaternionic Hessian and Monge-Amp\`ere operator. Now we are going
to describe them following the original approach of
\cite{alesker-bsm}. We also believe that in this language the
analogies with the classical real and complex cases become more
explicit.

We will write a quaternion $q\in \HH$ in the standard form
$$q= t+ x\cdot i +y\cdot j+ z\cdot k ,$$
where $t,\, x,\, y,\, z$ are real numbers, and $i,\, j,\, k$ satisfy
the usual quaternionic relations
$$i^2=j^2=k^2=-1, \, ij=-ji=k,\, jk=-kj=i, \, ki=-ik=j.$$

\def\db{\frac{\partial}{\partial \bar q}}
\def\dq{\frac{\partial}{\partial  q}}

The Dirac (or Cauchy-Riemann) operator $\frac {\partial}{\partial
\bar q}$ is defined as follows. For any $\HH$-valued function $F$
$$\db F:=\frac{\partial F}{\partial  t}  +
i \frac{\partial F}{\partial  x} + j \frac{\partial F}{\partial y} +
k \frac{\partial F}{\partial  z}.$$

Let us also define the operator $\dq$:
$$\dq F:=\overline{ \db \bar F}=
\frac{\partial F}{\partial  t}  -
 \frac{\partial F}{\partial x}  i-
 \frac{\partial F}{\partial  y} j-
\frac{\partial F}{\partial  z}  k.$$

In the case of several quaternionic variables, it is easy to see
that the operators $\frac{\pt}{\pt q_i}$ and $\frac{\pt}{\pt \bar
q_j}$ commute:
\begin{eqnarray}
\big[\frac{\pt}{\pt q_i},\frac{\pt}{\pt \bar q_j}\big]=0.
\end{eqnarray}

\hfill

For any real valued functions $f$ on the flat space $\HH^n$ the
matrix $\left(\frac{\pt^2 f}{\pt\bar q_i\pt q_j}\right)$ is
hyperhermitian; it corresponds exactly (up to a constant) to the
quaternionic Hessian. More precisely, using the isomorphism $t$ from
(\ref{E:isomor-t}) one has:
\begin{eqnarray}\label{E:hessian}
t(\pt\pt_J f)=\kappa\left(\frac{\pt^2 f}{\pt\bar q_i\pt 
q_j}\right),
\end{eqnarray}
where $\kappa>0$ is a normalizing constant, by Proposition 4.1 of
\cite{alesker-verbitsky-06}. The precise value of $\kappa$ will not
be important. In what follows it will be convenient to renormalize
the isomorphism $t$ to make this constant to be 1. We will denote by
$Hess_\HH f$ the matrix in the right hand side of (\ref{E:hessian})
(with $\kappa=1$).

It is not hard to show that a $C^2$ smooth function $f$ on $\HH^n$
is plurisubharmonic if and only if the hyperhermitian matrix
$\left(\frac{\pt^2 f}{\pt\bar q_i\pt q_j}\right)$ is non-negative
definite everywhere (see \cite{alesker-bsm},
\cite{alesker-verbitsky-06}).

\begin{proposition}[\cite{alesker-bsm}]
(i) Let $f:\HH ^n \to \HH$ be a smooth function. Then for any
$\HH$-linear transformation $A$ of $\HH ^n$ (as a right $\HH
$-vector space) one has the identities
$$ \left( \frac {\partial ^2 f(Aq)}{\partial \bar q_i \partial q_j} \right)
=A^* \left(\frac {\partial ^2 f}{\partial \bar q_i \partial q_j}(Aq)
\right)A .$$

(ii) If, in addition, $f$ is real valued then for any $\HH$-linear
transformation $A$ of $\HH ^n$ and any quaternion $a$ with $|a|=1$
$$ \left( \frac {\partial ^2 f(A(q \cdot a))}{\partial \bar q_i \partial q_j} \right)
=A^* \left(\frac {\partial ^2 f}{\partial \bar q_i \partial
q_j}(A(q\cdot a)) \right) A .$$
\end{proposition}

\hfill

It remains to rewrite the Monge-Amp\`ere operator $(\pt\pt_Jf)^n$ in
this language. Up to a positive normalizing constant which we
ignore, the Monge-Amp\`ere operator of a real valued function $f$ is
equal to the Moore determinant $\det \left(\frac{\pt^2 f}{\pt\bar q_i\pt
 q_j}\right)$.

\section{Second order estimate.}\label{S:second_order}
The main result of this section is Theorem \ref{T:second-order}
below. It establishes a uniform estimate on the Laplacian of the
solution of the Monge-Ampere equation (\ref{E:main-equation}). Let
us introduce a bit more notation. To shorten the notation, it will
be convenient to denote the quaternions $1,i,j,k$ by
$e_0,e_1,e_2,e_3$ respectively. Furthermore the $p$-th coordinate of
a quaternionic $n$-tuple $q=(q_1,\dots,q_n)$ will be written as
$$q_p=\sum_{i=0}^3 x^i_pe_p,$$
where $x^i_p\in \RR$. The partial derivative of a function $F$ with
respect to the real coordinate $x^i_p$ will be denoted by
$F_{x^i_p}$.

\hfill

First we prove the following elementary inequality.
\begin{proposition}\label{P:elem-inequality}
Let $u\in C^4$ be a strictly plurisubharmonic function such that at
a given point $z$ its quaternionic Hessian $(u_{\bar i j})$ is
diagonal. Then at this point $z$ one has
\begin{eqnarray}\label{C49}
\sum_{p=0}^3\sum_{i,k=1}^n\frac{|u_{\bar k kx_p^i}|^2}{u_{\bar i
i}u_{\bar k k}}\leq 2\sum_{p=0}^3\sum_{i,k,l=1}^n\frac{|u_{\bar k i
x_p^l}|^2}{u_{\bar i i}u_{\bar k k}}.
\end{eqnarray}
\end{proposition}
{\bf Proof.} Let us fix now the indices $i,k$ and compare the
summands containing this pair of indices in both sides of
(\ref{C49}).

First consider the case $i=k$. In the left hand side we have
\begin{eqnarray}\label{C50}
\sum_{p=0}^3\frac{|u_{\bar k kx_p^k}|^2}{u_{\bar k k}u_{\bar k k}}.
\end{eqnarray}
In the right hand side of (\ref{C49}) we have
\begin{eqnarray}\label{C51}
2\sum_{p=0}^3\sum_{l=1}^n\frac{|u_{\bar k kx_p^l}|^2}{u_{\bar k
k}u_{\bar k k}}
\end{eqnarray}
It is clear that $(\ref{C50})\leq (\ref{C51})$.

Let us consider the case now $i\ne k$. The left hand side of
(\ref{C49}) contains two summands with the pair $i,k$:
\begin{eqnarray}\label{C52}
\frac{1}{u_{\bar i i}u_{\bar k k}}\sum_p(|u_{\bar k
kx_p^i}|^2+|u_{\bar i ix_p^k}|^2).
\end{eqnarray}
The right hand side of (\ref{C49}) contains two summands with the
pair $i,k$:
\begin{eqnarray}\label{C53}
\frac{2}{u_{\bar i i}u_{\bar k k}}\sum_{p,l}(|u_{\bar k
ix_p^l}|^2+|u_{\bar i kx_p^l}|^2)=\frac{4}{u_{\bar i i}u_{\bar k
k}}\sum_{p,l}|u_{\bar k ix_p^l}|^2.
\end{eqnarray}
To finish the proof of proposition, it suffices now to show that
$(\ref{C52})\leq (\ref{C53})$, or explicitly after cancelling out
the term $u_{\bar i i}u_{\bar k k}$ on both sides, it reduces to
\begin{eqnarray}\label{C54}
\sum_p(|u_{\bar k kx_p^i}|^2+|u_{\bar i ix_p^k}|^2)\leq
4\sum_{p,l}|u_{\bar k ix_p^l}|^2.
\end{eqnarray}
In order to show such a general inequality it suffices to sum up in
the right hand side over $l=i,k$. Thus (\ref{C54}) follows from
\begin{eqnarray}\label{C55}
\sum_p(|u_{\bar k kx_p^i}|^2+|u_{\bar i ix_p^k}|^2)\leq
4\sum_p(|u_{\bar k ix_p^i}|^2+|u_{\bar k ix_p^k}|^2).
\end{eqnarray}
In the last inequality we may separate summands containing
derivatives $kki$ and $kii$. These two inequalities are completely
symmetric and obtained one from the other by exchange $i$ by $k$.
Thus it is enough to show
\begin{eqnarray}\label{C56}
\sum_p|u_{\bar k kx_p^i}|^2\leq 4\sum_p|u_{\bar k ix_p^k}|^2=4\sum_p|u_{\bar i kx^k_p}|^2.
\end{eqnarray}
Let us define two operators
$\overset{\leftarrow}{\pt_k},\overset{\leftarrow}{\pt_{\bar k}}$
acting on the space of quaternionic valued functions:
\begin{eqnarray*}
\overset{\leftarrow}{\pt_k}\Phi:=\sum_{p=0}^3\frac{\pt\Phi}{\pt
x_p^k}\bar e_p,\\
\overset{\leftarrow}{\pt_{\bar k}}\Phi:=\sum_{p=0}^3\frac{\pt
\Phi}{\pt x_p^k} e_p,
\end{eqnarray*}
where $\bar e_p$ denotes the quaternionic conjugate of the
quaternionic unit $e_p$ (here $p=0,...,3$).

Let $\Delta_k=\sum_{p=0}^3\frac{\pt^2}{(\pt x^k_p)^2}$ be the
Laplacian with respect to the $k$-th quaternionic variable. Clearly
$\Delta_k=\overset{\leftarrow}{\pt_k}\overset{\leftarrow}{\pt_{\bar
k}}=\overset{\leftarrow}{\pt_{\bar k}}\overset{\leftarrow}{\pt_k}$.
Also for real valued function $u$ $\Delta_ku=u_{\bar k k}$.

Let us take $\Phi=u_{\bar i}=\sum_qe_qu_{x_q^i}$. Then (\ref{C56})
is rewritten
\begin{eqnarray}\label{C57}
|\Delta_k\Phi|^2\leq 4
\sum_p|\overset{\leftarrow}{\pt_k}\Phi_{x_p^k}|^2
\end{eqnarray}
We have
$$\Delta_k\Phi=\overset{\leftarrow}{\pt_{\bar k}}\overset{\leftarrow}{\pt_k}\Phi=
\sum_q(\overset{\leftarrow}{\pt_k}\Phi)_{x_q^k} e_q.$$ Denote
$\Psi:=\overset{\leftarrow}{\pt_k}\Phi$. Then (\ref{C57}) is
rewritten:
\begin{eqnarray}\label{C58}
|\sum_q\Psi_{x_q^k} e_q|\leq 2\sqrt{\sum_q|\Psi_{x_q^k}|^2}
\end{eqnarray}

But indeed by the Cauchy-Schwarz $$|\sum_{q=0}^3\Psi_{x_q^k}
e_q|\leq \sum_{q=0}^3 |\Psi_{x_q^k}|\leq \sqrt{4}\cdot
\sqrt{\sum_{q=0}^3|\Psi_{x_q^k}|^2}.$$ \qed

\hfill

For any smooth function $g$ we denote by $g_a$, $g_{ab}$ the first
and second derivatives of $g$ with respect to coordinates with
indices $a$ and $a,b$ respectively (thus $a,b$ could be $x^i_p,
x^j_q$).
\begin{proposition}\label{D:det}
Let $U$ be a smooth function with values in hyperhermitian
invertible matrices. Let $ \det U= F.$

Then
\begin{eqnarray*}
Tr(U^{-1}U_{ab})=Tr(U^{-1}U_aU^{-1}U_b)+(\log F)_{ab}
\end{eqnarray*}

\end{proposition}
The proof is by straightforward computation using the identity
$$(\det U)_a=\det U\cdot Tr(U^{-1}U_a).$$

\hfill

We will need few more formulas. Let $U$ denote the quaternionic
Hessian of a function $u\in C^4$. Let $G$ be a smooth function with
values in positive definite hyperhermitian matrices. Define the
Laplacian
\begin{eqnarray}\label{D:A-Lap}
\Delta u:=Tr(G^{-1}U).
\end{eqnarray}
\begin{proposition}\label{P:dlaplacian}
Let $\gamma\colon \RR\to \RR$ be a smooth function. Then
\begin{eqnarray*}\label{E:dLaplacian}
[\gamma(\Delta u)]_{ab}=\\\gamma''(\Delta
u)\cdot[-Tr(G^{-1}G_aG^{-1}U)+Tr(G^{-1}U_a)]\cdot
[-Tr(G^{-1}G_bG^{-1}U)+Tr(G^{-1}U_b)]\\
+\gamma'(\Delta
u)\cdot\left[Tr(G^{-1}G_aG^{-1}G_bG^{-1}U)+Tr(G^{-1}G_bG^{-1}G_aG^{-1}U)\right.\\-
Tr(G^{-1}G_{ab}G^{-1}U)\\-Tr(G^{-1}G_aG^{-1}U_b)-Tr(G^{-1}G_bG^{-1}U_a)\\\left.+Tr(G^{-1}U_{ab})\right]
\end{eqnarray*}
\end{proposition}
{\bf Proof.} We have
\begin{eqnarray*}
[\gamma(\Delta u)]_{ab}=[\gamma'(\Delta u)\cdot (\Delta u)_a]_b=\\
\gamma''(\Delta u) \cdot (\Delta u)_a \cdot (\Delta
u)_b+\gamma'(\Delta u)\cdot (\Delta u)_{ab}.
\end{eqnarray*}
Next we have
\begin{eqnarray*}
(\Delta u)_a=(Tr(G^{-1} U))_a=\\
-Tr(G^{-1}G_aG^{-1} U)+Tr(G^{-1}U_a).
\end{eqnarray*}
Also
\begin{eqnarray*}
(\Delta u)_{ab}=[-Tr(G^{-1}G_aG^{-1} U)+Tr(G^{-1}U_a)]_b=\\
Tr(G^{-1}G_aG^{-1}G_bG^{-1}U)+Tr(G^{-1}G_bG^{-1}G_aG^{-1}U)\\-
Tr(G^{-1}G_{ab}G^{-1}U)\\-Tr(G^{-1}G_aG^{-1}U_b)-Tr(G^{-1}G_bG^{-1}U_a)\\+Tr(G^{-1}U_{ab}).
\end{eqnarray*}
The proposition follows. \qed

\hfill

Given a fixed plurisubharmonic function $u\in C^4$, let us define
another Laplacian
\begin{eqnarray}\label{D:B-Lap}
\Delta ' v:=Tr(U^{-1}V)
\end{eqnarray}
where $U$ and $V$ are quaternionic Hessians of $u$ and $v$
respectively. Proposition \ref{P:dlaplacian} implies immediately
\begin{proposition}\label{P:flat-double-Laplacian}
Let $u\in C^4$ be a strictly plurisubharmonic function. Let us
assume that $G$ is a flat hyperK\"ahler metric. Then choose
(locally) coordinates such that $G\equiv Id$. Then in these
coordinates
\begin{eqnarray*}
\,\,[\gamma(\Delta u)]_{ab}=\gamma''(\Delta u)\cdot Tr(U_a)\cdot
Tr(U_b)+\gamma'(\Delta u)Tr(U_{ab}).
\end{eqnarray*}
If moreover at a point $z$ the matrix $U(z)=(u_{\bar i j}(z))$ is
diagonal then at this point $z$ one has
\begin{eqnarray*}
\Delta'(\gamma(\Delta u))=\gamma''(\Delta
u)\sum_{i,p}\frac{1}{u_{\bar i i}}(TrU_{x_p^i})^2+\gamma'(\Delta
u)\sum_{i,k}\frac{u_{\bar i i \bar k k}}{u_{\bar i i}}.
\end{eqnarray*}
\end{proposition}

\begin{corollary}\label{C:Lapl-Lapl}
Let $u\in C^4$ be a strictly plurisubharmonic function. Let us
assume that $G$ is a flat metric. Let us fix a point $z$. Then
choose (locally) coordinates such that $G\equiv Id$ in a
neighborhood, and $U:=(u_{\bar i j})$ is diagonal at $z$. Let
$F:=\det U$ as previously. Then in these coordinates we have at the
point $z$
\begin{eqnarray*}
\Delta'(\gamma(\Delta u))=\\\gamma''(\Delta
u)\sum_{i,p}\frac{1}{u_{\bar i i}}(TrU_{x_p^i})^2 +\gamma'(\Delta
u)\left[\sum_{i,p}Tr((U^{-1}U_{x^i_p})^2) +\Delta(\log F)\right]=\\
\gamma''(\Delta u)\sum_{i,p}\frac{1}{u_{\bar i i}}(\sum_k u_{\bar k k
x^i_p})^2 +\gamma'(\Delta u)\left[\sum_{i,l,n,p}\frac{|u_{\bar l n
x^i_p}|^2}{u_{\bar l l}u_{\bar n n}}+\sum_i(\log F)_{\bar i i}\right],
\end{eqnarray*}
where $\Delta$ and $\Delta'$ are defined by (\ref{D:A-Lap}) and
(\ref{D:B-Lap}) respectively.
\end{corollary}
{\bf Proof.} The second equality is just immediate substitution of
matrix $U=diag(u_{\bar 1 1},\dots, u_{\bar n n})$. Let us prove the
first one. By Proposition \ref{P:flat-double-Laplacian} it suffices
to show that
\begin{eqnarray*}
\sum_{i,k}\frac{u_{\bar i i \bar k k}}{u_{\bar i
i}}=\sum_{i,p}Tr((U^{-1}U_{x^i_p})^2)+\Delta(\log F).
\end{eqnarray*}
The left hand side of the last equality is equal to $\sum_{k,p}
Tr(U^{-1}U_{x^k_px^k_p}).$ But by Proposition \ref{D:det} the last
expression is equal to
\begin{eqnarray*}
\sum_{k,p}\left(Tr(U^{-1}U_{x^k_p}U^{-1}U_{x^k_p})+(\log
F)_{x^k_px^k_p}\right)=\sum_{k,p}Tr((U^{-1}U_{x^k_p})^2)+\sum_k(\log
F)_{\bar k k}.
\end{eqnarray*}
\qed

\begin{proposition}\label{P:ineq-Lapl-psh}
Let $u\in C^4$ be a strictly plurisubharmonic function. Denote $\det
U=F$ as previously. Let $G$ be a locally flat quaternionic metric.
Then
$$\Delta'(2\sqrt{\Delta u})\geq (\Delta u)^{-1/2}\Delta(\log F),$$
where $\Delta$ and $\Delta'$ are defined by (\ref{D:A-Lap}) and
(\ref{D:B-Lap}) respectively.
\end{proposition}
{\bf Proof.} We prove it pointwise. Let us fix a point $z$. We can
choose coordinates near $z$ such that $G\equiv Id$ and $U(z)$ is
diagonal. Then clearly $\Delta h=\sum_i h_{\bar i i}$.

Let us take in Corollary \ref{C:Lapl-Lapl} $\gamma(x)=2\sqrt{x}$.
Then
\begin{eqnarray*}
\Delta'(2\sqrt{\Delta u})=(\Delta u)^{-1/2}\sum_i(\log F)_{\bar i
i}\\+(\Delta u)^{-1/2}\left[\sum_{i,l,n,p}\frac{|u_{\bar l n
x^i_p}|^2}{u_{\bar l l}u_{\bar n n}}-\frac{1}{2\Delta
u}\sum_{i,p}\frac{1}{u_{\bar i i}}(\sum_k u_{\bar k k
x^i_p})^2\right].
\end{eqnarray*}
It remains to show that the expression in the square brackets in
non-negative. We are using the Cauchy-Schwarz inequality
\begin{eqnarray*}
\frac{1}{2\Delta u}\sum_{i,p}\frac{1}{u_{\bar i i}}(\sum_k u_{\bar k k
x^i_p})^2=\frac{1}{2\Delta u}\sum_{i,p}\frac{1}{u_{\bar i
i}}(\sum_k\sqrt{u_{\bar k k}}\frac{u_{\bar k k x^i_p}}{\sqrt{u_{\bar k
k}}})^2\leq\\
\frac{1}{2}\sum_{i,p}\frac{1}{u_{\bar i i}}\sum_k\frac{|u_{\bar k k
x^i_p}|^2}{u_{\bar k k}}=\frac{1}{2}\sum_{i,k,p}\frac{|u_{\bar k k
x^i_p}|^2}{u_{\bar i i}u_{\bar k k}}.
\end{eqnarray*}
But by Proposition \ref{P:elem-inequality} the last expression does
not exceed $\sum_{i,l,n,p}\frac{|u_{\bar l n x^i_p}|^2}{u_{\bar l
l}u_{\bar n n}}$. \qed

\hfill

Now we return back to the Monge-Amp\`ere equation.
\begin{theorem}\label{T:second-order}
Let $M$ be a compact manifold with a locally flat hypercomplex
structure. Let us assume in addition that $M$ admits a metric $G$
which is parallel with respect to the Obata connection.
\footnote{Any such metric $G$ parallel with respect to the Obata
connection is automatically hyperK\"ahler. Hence equivalently one
can say that $M$ admits a locally flat hyperK\"ahler metric
compatible with the hypercomplex structure.} Let $G_0$ be another
HKT-metric on $M$. Let $\phi\colon M\to \RR$ be a $C^4$-smooth
solution of the quaternionic Monge-Amp\`ere equation
\begin{eqnarray}\label{E:monge-ampere}
\det(G_0+Hess_\HH\phi)=e^f\det G_0
\end{eqnarray}
where $f$ is a $C^2$-smooth function. Then there exists a constant
$C$ depending on $M,G,G_0$, and $||f||_{C^2(M)}$ such that
$$||\Delta_G\phi||_{C^0}\leq C$$
where $\Delta_G\colon C^2(M,\RR)\to \RR$ is the globally defined
operator which in flat local coordinates is equal $\Delta _G
h:=Tr(G^{-1}\cdot Hess_{\HH}(h))$.
\end{theorem}
{\bf Proof.} Let us denote
$$\Delta'h:= Tr((G_0+Hess_\HH\phi)^{-1}\cdot Hess_\HH(h)).$$

Let $\Ome$ and $\Ome_0$ be the HKT-forms corresponding to $G$ and
$G_0$ respectively. We may assume that the solution $\phi$ satisfies
$$\int_M \phi\cdot \Ome_0^n\wedge\bar \Ome^n=0.$$
Then by Corollary 5.7 in \cite{alesker-verbitsky-10} (the uniform
estimate), there exists a constant $C_1$ such that
$||\phi||_{C^0}\leq C_1$. Let us consider the function
$$T:=2\sqrt{Tr(G^{-1}\cdot (G_0+Hess_\HH(\phi)))}-\phi.$$ In order to
prove the theorem, it suffices to show that this function is bounded
from above. Let $z\in M$ be a point of maximum of the function $T$.
Then
\begin{eqnarray}\label{E:01}
\Delta'T(z)\leq 0.
\end{eqnarray}
Let $g_0\in C^\infty$ be a local potential of the metric $G_0$. Then
$u:=g_0+\phi\in C^4$ is a strictly plurisubharmonic function. Let
$U$ denote its quaternionic Hessian. In flat local coordinates
around $z$ we can rewrite the Monge-Amp\`ere equation
(\ref{E:monge-ampere}) as
\begin{eqnarray}\label{E:MA2}
\det U=F
\end{eqnarray}
where $F$ is identified with $e^f\det G_0$. Also in this notation
$T=2\sqrt{\Delta_G u}-\phi$. From this and (\ref{E:01}) we get
\begin{eqnarray}\label{E:02}
(\Delta'(2\sqrt{\Delta_G u}-\phi))(z)\leq 0.
\end{eqnarray}
By (\ref{E:02}) and Proposition \ref{P:ineq-Lapl-psh} we obtain
\begin{eqnarray}
0\geq (\Delta_G u)^{-1/2}\sum_i(\log F)_{\bar i i}-
\Delta'\phi=\\
(\Delta_G u)^{-1/2}\sum_i(\log F)_{\bar i i}- Tr(U^{-1}\cdot
(U-G_0))=\\
(\Delta_G u)^{-1/2}\sum_i(\log F)_{\bar i i}-n+ Tr(U^{-1}\cdot G_0).
\end{eqnarray}
Let us choose coordinates near $z$ so that $G\equiv Id$ and $U(z)$
is diagonal. Let $C_2:=n||\log F||_{C^2}$. Then we get
\begin{eqnarray}\label{E:03}
0\geq -C_2\cdot (\Delta_G u)^{-1/2}-n+ Tr(U^{-1}\cdot
G_0)\\\label{E:04} -C_2\cdot (\sum_i u_{\bar i
i}(z))^{-1/2}-n+\sum_i\frac{(g_0)_{\bar i i}(z)}{u_{\bar i i}(z)}.
\end{eqnarray}
Let $C_3>0$ be a constant, depending on $M$ and $G_0$ only, such
that for all $i=1,\dots, n$ one has
$$(g_0)_{\bar i i}(z)\geq C_3.$$ This and (\ref{E:03})-(\ref{E:04})
imply
\begin{eqnarray}\label{E:05}
C_3\cdot \sum_i\frac{1}{u_{\bar i i}(z)}\leq C_2\cdot (\sum_i
u_{\bar i i}(z))^{-1/2}+n.
\end{eqnarray}
Since $\prod_i u_{\bar i i}=F$, the arithmetic-geometric mean
inequality implies that
\begin{eqnarray}\label{E:06}
(\sum_i u_{\bar i i})^{-1/2}\leq n^{-1/2}F^{-1/2n}.
\end{eqnarray}
Let $C_4=C_2\cdot n^{-1/2}||F^{-1/2n}||_{C^0}$. Hence we get
\begin{eqnarray}\label{E:07}
C_3\cdot \sum_i\frac{1}{u_{\bar i i}(z)}\leq C_4 +n.
\end{eqnarray}
Obviously (\ref{E:07}) implies that $\sum_i\frac{1}{u_{\bar i
i}(z)}\leq C_5$, and hence $$u_{\bar i i}(z)>C_6>0 \mbox{ for all }
i.$$ Since $\prod_i u_{\bar i i}=F$ we obtain
$$\sum_i u_{\bar i i}(z)\leq C_7.$$ Hence $(\Delta_G u)(z)\leq C_7$.
This implies the proposition. \qed

\begin{remark}
If $u\in C^2(M)$ then estimates on $||u||_{C^0(M)}$ and $||\Delta_G
u||_{C^0(M)}$, imply an estimate on $||u||_{C^{1,\alp}(M)}$ for any
$0<\alp<1$ by Theorem 8.32 in \cite{gilbarg-trudinger}.
\end{remark}

\begin{remark}
If we replace $G$ by any HKT-metric $G_1$, then we can define
similarly the operator $\Delta_{G_1}h:= Tr(G_1^{-1}\cdot Hess_\HH
h)$. By a simple linear algebra it is easy to show that an estimate
on $||\Delta_Gu||_{C^0(M)}$ is equivalent to an estimate on
$||\Delta_{G_1}u||_{C^0(M)}$
\end{remark}

\section{$C^{2,\alp}$-estimate.}\label{S:2-alpha}
Let $(M^{4n}, I,J,K)$ be a compact hypercomplex manifold. Let
$\Ome_0\in \Lam^{2,0}_I$ be an HKT-form.
We are interested in the quaternionic Monge-Amp\`ere equation
\begin{eqnarray}\label{MA-main}
(\Ome_0+\pt\pt_J\phi)^n=e^f\Ome_0^n.
\end{eqnarray}

Let us denote
$$\Delta\phi=\frac{\pt\pt_J\phi\wedge \Ome_0^{n-1}}{\Ome_0^n}.$$
Clearly $\Delta$ is a linear second order elliptic operator without
free term (i.e. $\Delta(1)=0$) and with infinitely smooth
coefficients.

The main result of this section is as follows.
\begin{theorem}\label{thm-main}
Let $M^{4n}$ be a compact manifold with locally flat hypercomplex
structure. There exist $\alp\in (0,1)$ and $C>0$, both depending on
$M$,$\Ome_0$,
$||f||_{C^2(M)}$,$||\phi||_{C^0(M)}$,$||\Delta\phi||_{C^0(M)}$ only,
such that
$$||\phi||_{C^{2,\alp}(M)}\leq C.$$
\end{theorem}

Recall that by Theorem \ref{m}(2) locally $\Ome_0$ can be
represented by a potential $\Ome_0=\pt\pt_J g_0$ where $g_0\in
C^\infty$ is quaternionic strictly plurisubharmonic. Since $M$ is
locally flat, Theorem \ref{thm-main} follows from the following
version on $\HH^n$ applied to $u=g_0+\phi$.
\begin{theorem}\label{thm-local}
Let $u\in C^4$ be a quaternionic psh function in an open subset
$\co\subset \HH^n$ satisfying
\begin{eqnarray}\label{2}
\det(u_{\bar i j})=e^f
\end{eqnarray}
with $f\in C^2$. Let $\co'\subset \co$ be a relatively compact open
subset. Then there exist $\alp\in(0,1)$ depending on
$n,||u||_{C^0(\co)},||\Delta u||_{C^0(\co)},$ $ ||f||_{C^2(\co)}$
only and a constant $C>0$ depending in addition on
$dist(\co',\pt\co)$ such that
$$||u||_{C^{2,\alp}(\co')}\leq C.$$
\end{theorem}
The proof of this theorem is a quaternionic version of the
Evans-Krylov method \cite{evans}-\cite{evans2}, \cite{krylov}. The
complex version of it was considered by Siu \cite{siu} and B{\l}ocki
\cite{blocki-duke-2000}. Our exposition closely follows B{\l}ocki
\cite{blocki-lecture-notes}; perhaps only Lemma \ref{L:divergence}
below is somewhat novel.

\hfill

\def\dz{\Delta_\zeta}

For a unit vector $\zeta\in \HH^n$ we denote by $\dz$ the Laplacian
on any translate of the (right) quaternionic line spanned by
$\zeta$. Also let us denote by $U$ the quaternionic Hessian
$(u_{\bar i j})$. Thus $U$ is a hyperhermitian positive definite
$n\times n$ matrix.

\begin{lemma}\label{L:p-kvadrat}
Assume that $u,f$ satisfy the assumptions of Theorem
\ref{thm-local}. Then pointwise we have
\begin{eqnarray}
Tr(U^{-1}\cdot\dz U)\geq \dz f.
\end{eqnarray}
\end{lemma}
\def\xp1{x_p^1x_p^1}
{\bf Proof.} We may assume that $\zeta=(1,0,\dots,0)$. Then
$\dz=\sum_{p=0}^3\frac{\pt^2}{(\pt x_p^1)^2}$. It is enough to show
that
\begin{eqnarray}\label{4}
Tr(U^{-1}\cdot U_{x_p^1 x_p^1})\geq f_{x_p^1 x_p^1}
\end{eqnarray}
for any $p=0,\dots, 3$. Differentiating the equality
$$\log\det U=f$$
twice with respect to $x_p^1$, we obtain
\begin{eqnarray*}\label{5}
Tr(U^{-1}\cdot U_{\xp1})=
f_{\xp1}+Tr(U^{-1}U_{x_p^1}U^{-1}U_{x_p^1}).
\end{eqnarray*}
In order to prove the lemma, it suffices to show that
$Tr(U^{-1}U_{x_p^1}U^{-1}U_{x_p^1})\geq 0$. More generally let us
show that if $A,B$ are hyperhermitian matrices and $A>0$ then
$$Tr(A^{-1}BA^{-1}B)\geq 0.$$
Since $A>0$ we can diagonalize $A,B$ simultaneously. More precisely
we can find an invertible quaternionic matrix $T$ and a real
diagonal matrix $D$ such that
$$A=T^*T,\, B=T^*DT.$$
Then
\begin{eqnarray*}
Tr(A^{-1}BA^{-1}B)=Tr(T^{-1}D^2 T)=Tr(D^2)\geq 0.
\end{eqnarray*}
Lemma is proved. \qed

\hfill

Let us recall now the weak Harnack inequality (see
\cite{gilbarg-trudinger}, Theorem 8.18, or \cite{han-lin} Theorem
4.15). Below we normalize everywhere the Lebesgue measure on $\RR^N$
by $vol(B_1)=1$ where $B_1$ is Euclidean ball of unit radius. We
also denote $D_i=\frac{\pt}{\pt x_i}$.

\begin{theorem}[weak Harnack inequality]\label{harnack}
Let $B_R\subset \RR^N$ be a Euclidean ball of radius $R$. Let
$(a_{ij})_{i,j=1}^N\in L^\infty(B_R)\cap C^1(B_R)$, $a_{ij}=a_{ji}$,
satisfy uniform elliptic estimate
$$\lam||\xi||^2\leq \sum_{i,j}a_{ij}(x)\xi_i\xi_j\leq
\Lam||\xi||^2,\mbox{ for all }\xi\in \RR^N$$ with $\lam,\Lam>0$. Let
$v\in C^2(B_R)$ be a function satisfying
\begin{eqnarray*}
v\geq 0,\\
\sum_{i,j}D_j(a_{ij}D_iv)\leq \psi,
\end{eqnarray*}
where $\psi\in L^\infty(B_R)$. Then for any $0<\theta<\tau<1$ we
have
$$\inf_{B_{\theta R}}v+R||\psi||_{L^\infty(B_R)}\geq
C\left(R^{-N}\int_{B_{\tau R}}v\right)$$ where the constant $C$
depends only on $\lam,\Lam,\theta,\tau,N$.
\end{theorem}
\begin{remark}\label{R:harnack}
We will use Theorem \ref{harnack} in the following weaker form. We
will take $R=4r, \theta=1/4,\tau=1/2$. Then we deduce
\begin{eqnarray}
r^{-N}\int_{B_r}v\leq C'(\inf_{B_r} v+r)
\end{eqnarray}
where the constant $C'$ depends on
$\lam,\Lam,||\psi||_{L^\infty(B_R)},N$ only.
\end{remark}

\hfill

For $U=(u_{\bar i j})$ as above define the operator $\cd$ by
\begin{eqnarray}\label{D:D-oper}
\cd h=\det U\cdot Tr(U^{-1}\cdot Hess_\HH h).
\end{eqnarray}

First let us prove an algebraic lemma.
\begin{lemma}\label{L:divergence}
The operator $\cd$ defined by (\ref{D:D-oper}) can be written in the
divergence form as in Theorem \ref{harnack}, namely
$$\cd h=\sum_{st}D_s(a_{st}D_th)$$
where $a_{st}$ is a $4n\times 4n$ real symmetric matrix with
$C^2$-smooth coefficients, and $s,t$ in the sum run over all real
variables $x^i_p$.
\end{lemma}
Before we prove the lemma, let us prove the following linear
algebraic claim.
\begin{claim}\label{Cl-LinAlg}
Let $A,B$ be $n\times n$ hyperhermitian matrices. Suppose that $A$
is invertible. Then
\begin{eqnarray}\label{E:LinAlg}
\det A\cdot Tr(A^{-1} B)=n \det(A[n-1],B).
\end{eqnarray}
\end{claim}
{\bf Proof.} Both sides of the equality are linear in $B$. Hence it
suffices to prove the equality for $B>0$. Then $A$ and $B$ can be
diagonalized simultaneously, more precisely the exists an invertible
quaternionic matrix $T$ and a real diagonal matrix $D$ such that
$$B=T^*T,\, A=T^*DT.$$
Then the left hand side of (\ref{E:LinAlg})  is equal to
\begin{eqnarray*}
\det(T^*DT)\cdot Tr(T^{-1}D^{-1}T)=\det(T^*T)\cdot\det D\cdot
Tr(D^{-1}).
\end{eqnarray*}
On the other hand the right hand side is equal to
$$n\det((T^*DT)[n-1],T^*T)=n\det(T^*T)\cdot
\det(D[n-1],I_n).$$ Hence it suffices to assume that $B=I_n$ and
$A=D$ is real diagonal, i.e.
$$\det D\cdot Tr(D^{-1})=n\det(D[n-1],I_n).$$
The last identity for real diagonal $D$ is obvious. \qed

\hfill

{\bf Proof of Lemma \ref{L:divergence}.} Let us consider on $\HH^n$
the complex structure $I$.

By \cite{alesker-verbitsky-06}, Corollary 4.6, for appropriate
choice of flat $I$-complex coordinates on $\HH^n$ one has
$$(\pt\pt_J h)^n= \kappa_n \det (h_{\bar i j}) \cdot dz_1\wedge dz_2\wedge\dots\wedge dz_{2n}$$
where $\kappa_n$ is a normalizing constant depending on $n$ only
(its precise value will not be important in the argument below).
Polarizing the last equality we obtain for any $n$-tuple of
functions $h_1,\dots,h_n$
\begin{eqnarray}\label{E:1-1}
(\pt\pt_J h_1)\wedge\dots \wedge (\pt\pt_J h_n)=\det((h_1)_{\bar i
j},\dots,(h_n)_{\bar i j})\Theta,
\end{eqnarray}
where we have introduced the notation $\Theta:=\kappa_n dz_1\wedge
dz_2\wedge\dots\wedge dz_{2n}\in \Lam^{2n,0}_{I}(\HH^n)$ for
brevity. Hence
\begin{eqnarray}
\cd h=\det U\cdot Tr(U^{-1}(h_{\bar i j}))\overset{\mbox{Claim }
\ref{Cl-LinAlg}}{=}\\
n\det(U[n-1],(h_{\bar i
j}))\overset{(\ref{E:1-1})}{=}\\n\label{E:1-2} \frac{(\pt\pt_J
u)^{n-1}\wedge \pt\pt_J h}{\Theta}.
\end{eqnarray}

Let $(a_{st})$ be the $4n\times 4n$ be the real symmetric matrix
defined to be the realization of the $n\times n$ quaternionic
hermitian matrix $\det U\cdot U^{-1}$. Then it is easy to see that
$\cd h=\sum_{s,t} a_{st}D_sD_t h$.

Clearly the statement of Lemma \ref{L:divergence} is equivalent to
\begin{eqnarray}\label{E:1-3}
\sum_s D_s a_{st}=0 \mbox{ for any }t.
\end{eqnarray}
In order to prove the last equality let us rewrite it in a more
invariant way. Let $\nabla$ denote the flat connection on the
tangent bundle of $\HH^n=:M$. Let $q\colon T^*M\otimes TM\to \RR$ be
the natural pairing. Let
$$Q\colon Sym^2TM\otimes T^*M\to TM$$
be the natural contraction map given by $Q(x\otimes y\otimes
\xi)=\xi(y) x$.

It is clear that the quaternionic Hessian $U=(u_{\bar i j})$ belongs
to the space $\cb$ of quadratic forms on $\HH^n$ which are invariant
under the (right) multiplication by norm one quaternions. Hence the
matrix $a:=(a_{st})$, which corresponds to $\det U\cdot U^{-1}$,
belongs to $\cc:=\cb^*\otimes L$ where $L$ denotes the line to which
the Moore determinant belongs (below we will identify $L$ more
explicitly).

In this notation (\ref{E:1-3}) is equivalent to
\begin{eqnarray}\label{E:1-4}
Q(\nabla a)=0.
\end{eqnarray}
Since $a$ changes as an appropriate tensor under all translations on
$\HH^n$ and all linear transformations from $GL_n(\HH)\cdot
GL_1(\HH)$, and since $\nabla$ commutes with such transformations,
the equation (\ref{E:1-4}) is invariant under the group
$\HH^n\rtimes (GL_n(\HH)\cdot GL_1(\HH))$. Hence it suffices to
check it at the point 0. This is the first order differential
equation. The 1-jet of $a$ at 0 belongs to the space $\cc\bigoplus
\cc\otimes_\RR (\HH^n)^*$. The differential operator $a\mapsto
Q(\nabla a)$ obviously does not depend on the first component of $a$
in this sum. Thus let us denote by $j(a)$ the second component of
$a$. The subspace of $\cc\otimes_\RR (\HH^n)^*$ corresponding to
solutions of the equation (\ref{E:1-4}) is a $GL_n(\HH)\cdot
GL_1(\HH)$-invariant proper subspace (clearly it is not equal to the
whole space, and that it is non-zero will be seen from the last part
of the argument where we will construct non-zero examples of
solutions of this equation). Let us study the decomposition of
$\cc\otimes_\RR (\HH^n)^*$ under the group $GL_n(\HH)\cdot
GL_1(\HH)$. Actually it will be convenient to replace this group by
$GL_n(\HH)\times GL_1(\HH)$ which is mapped surjectively onto it.
Also it will be convenient to replace all spaces and groups by their
complexifications. We have
$$\HH^n\otimes_\RR \CC=V\otimes_\CC W$$
where $V=\CC^{2n},\, W=\CC^2$. It is well known (and easy to see
directly) that
$$\cb\otimes_\RR \CC=Sym^2V^*\otimes \det W^*,$$
where in the right hand side $Sym$, $\wedge$, and $\otimes$ are
taken over $\CC$ (here and below we will omit this subscript). It is
easy to see that the complexified line $L$ where the Moore
determinant lies is equal to $L\otimes_\RR\CC=\det V^*\otimes (\det
W^*)^{\otimes n}$. Hence we obtain
\begin{eqnarray}\label{E:1-5}
\left(\cc\otimes_\RR (\HH^n)^*\right)\otimes_\RR \CC=\\(Sym^2
V\otimes \det W)\otimes V^*\otimes W^*\otimes \left(\det V^*\otimes
(\det W^*)^{\otimes n}\right) \end{eqnarray} where all the tensor
products on the right hand side are over $\CC$.

Next the complexification of the group $GL_n(\HH)\times GL_1(\HH)$
is equal to the product  $GL(V)\times GL(W)$ of complex linear
groups, where $GL(V)$ acts on $V$ in the standard way, and similarly
for $W$.

Obviously $\det W\otimes W^*$ is an irreducible
$GL(W)$-representation. However the $GL(V)$-representation
$Sym^2V\otimes V^*$ is a direct sum of two irreducible
non-isomorphic subspaces (this easily follows from the Schur-Weyl
duality, see e.g. \cite{goodman-wallach}, \S 9.1.1). Hence the
representation of $GL(V)\times GL(W)$ in the space (\ref{E:1-5}) is
also a direct sum of two irreducible non-isomorphic subspaces.

The complexification of the map $Q$ tensored with the identity map
of $L\otimes_\RR\CC$, which we also will denote by $Q$, maps
\begin{eqnarray*}
Q\colon (Sym^2 V\otimes \det W)\otimes V^*\otimes W^*\otimes
\left(\det V^*\otimes (\det W^*)^{\otimes n}\right)\to \\V\otimes
W\otimes\left(\det V^*\otimes (\det W^*)^{\otimes n}\right).
\end{eqnarray*}
It is equal to the tensor product of the obvious contraction maps
\begin{eqnarray*}
Sym^2V\otimes V^*\to V,\\
\det W\otimes W^*=\wedge^2 W\otimes W^*\to W.
\end{eqnarray*}
Hence the kernel of $Q$ is an irreducible non-zero $GL(V)\times
GL(W)$-representation which we will denote by $\ck$.

\hfill

Thus to prove the lemma it remains to show that $j(a)\in \ck$.
Recall that $a$ corresponds to the hyperhermitian matrix
$\det(Hess_\HH (u))\cdot Hess_\HH(u)$ where $u$ is a function. Due
to (\ref{E:1-2}) and some linear algebra, this expression can be
identified with $(\pt\pt_J u)^{n-1}$ (up to a constant). But the
last expression satisfies
$$\pt((\pt\pt_J u)^{n-1})=\pt_J((\pt\pt_J u)^{n-1})=0.$$
These equations give a non-trivial restriction on $j(a)$ and imply
that $j(a)$ belongs to a proper $GL(V)\times GL(W)$-invariant
subspace for all functions $u$. Hence $j(a)$ belongs either always
(i.e. for any $a$) to $\ck$ or always to the other irreducible
summand of the space (\ref{E:1-5}). Let us show that the first case
takes place.

\hfill

It suffices to give an example of a function $u$ such that the
corresponding $a$ satisfies $Q(\nabla a)=0$ and $j(a)\ne 0$. Let us
write
$$\HH^n=\RR^n\oplus \RR^n\cdot i\oplus \RR^n\cdot j\oplus \RR^n\cdot
k.$$ We take an arbitrary smooth strictly convex function $u$ on the first copy on $\RR^n$
and extend it to the whole $\HH^n$ so that the extension is independent of the other $3n$ coordinates.
Then the quaternionic
Hessian is $u$ is equal to its usual real Hessian with respect to
the first $n$ real coordinates:
$$Hess_\HH u=\left(\frac{\pt^2u}{\pt t_\alp\pt
t_\beta}\right)=:(u_{\alp\beta})_{n\times n}.$$ Then set
$$A=\det (u_{\alp\beta})\cdot (u_{\alp\beta})^{-1}.$$
Then by definition $a$ is equal to the realization of $A$. Hence
\begin{eqnarray*}
a=\left[\begin{array}{cccc}
         A&0&0&0\\
         0&A&0&0\\
         0&0&A&0\\
         0&0&0&A
        \end{array}\right]_{4n\times 4n}.
\end{eqnarray*}
Then our equation $Q(\nabla a)=0$ is equivalent to
$$\frac{\pt}{\pt t_\alp}a_{\alp\beta}=0,$$
where we have summation with respect to repeated indices. More
explicitly  the last equation can be rewritten
\begin{eqnarray}\label{E:1-6}
\frac{\pt}{\pt t_\alp}(\det U\cdot U^{\alp\beta})=0
\end{eqnarray}
where $U=(u_{\alp\beta}),\, U^{\alp\beta}=(U^{-1})_{\alp\beta}$.
Next we have
\begin{eqnarray*}
\frac{\pt}{\pt t_\alp}\det U=\det U \cdot Tr(U^{-1} \frac{\pt U}{\pt
t_\alp}),\\
\frac{\pt U^{-1}}{\pt t_\alp}=-U^{-1}\frac{\pt U}{\pt t_\alp}U^{-1}.
\end{eqnarray*}
Hence, again with a summation over repeated indices, we get
\begin{eqnarray*}
\frac{\pt}{\pt t_\alp}(\det U\cdot U^{\alp\beta})=\det
U\left(U^{pq}u_{qp\alp}U^{\alp\beta}-U^{\alp
p}u_{pq\alp}U^{q\beta}\right)=0.
\end{eqnarray*}
Hence $j(a)\in \ck$ for any $u$ as above. It is easy to see that $u$
can be chosen so that $j(a)\ne 0$. This implies the lemma. \qed


\hfill

Next let us observe that $\cd$ is uniformly elliptic with constants
$\lam,\Lam$ depending on $||f||_{C^0(\co)},
||\Delta\phi||_{C^0(\co)}$ only. We are going to apply the weak
Harnack inequality to the operator $\cd$ to the function
$$v:=\sup_{B_{4r}}\dz u-\dz u$$
where $B_{4r}=B(z_0,4r)\subset \co$ is a ball with a center $z_0\in
\co'$.

The function $v$ satisfies by Lemma \ref{L:p-kvadrat}
$$\cd v\leq -e^{f}\dz f$$
where we have used $\Delta_\zeta(u_{\bar i j})=(\Delta_\zeta
u)_{\bar i j}$.

Hence we can apply the weak Harnack inequality in the form given in
Remark \ref{R:harnack}:
\begin{eqnarray}\label{harnack-main}
r^{-4n}\int_{B_r}(\sup_{B_{4r}}\dz u-\dz u)\leq C(\sup_{B_{4r}}\dz
u-\sup_{B_r}\dz u +r)
\end{eqnarray}
where $C$ depends on $||f||_{C^2(\co)},||\Delta \phi||_{C^0(\co)},
n$ only.

\begin{lemma}\label{L:det-ineq}
Let $A,B$ be positive definite hyperhermitian matrices of size $n$.
Then
$$Tr(A^{-1}(A-B))\leq n(\det A)^{-1/n}((\det A)^{1/n}-(\det B)^{1/n}).$$
\end{lemma}
{\bf Proof.} $A$ and $B$ can be simultaneously diagonalized. Thus
the inequality follows from the corresponding result in the real
diagonal case when it is equivalent to the arithmetic-geometric mean
inequality. \qed

\hfill

Lemma \ref{L:det-ineq} implies that for any point $x,y\in \co$
\begin{eqnarray}\label{E:II}
e^{f(y)}Tr\left(U^{-1}(y)(U(y)-U(x))\right)\leq \\\label{E:II-1}
ne^{f(y)\cdot (1-\frac{1}{n})}(e^{f(y)/n}-e^{f(x)/n})\leq
C_1||x-y||,
\end{eqnarray}
where $C_1$ depends on $||f||_{C^1(\co)}$ and $n$.

\hfill

Let $\zeta\in \HH^n$ be a vector. If we write it as a column, it
defines a rank one hyperhermitian matrix $\zeta\otimes \zeta^*$.
Namely $(\zeta\otimes\zeta^*)_{i\bar \j}=\zeta_i\cdot \bar \zeta_j$.
Clearly the matrix $\zeta\otimes \zeta^*$ is non-negative definite.
We will need the following linear algebraic lemma which is
completely analogous to the real and complex cases (see
\cite{gilbarg-trudinger}, Lemma 17.13, for the real case, and for
the complex case \cite{siu}, p. 103, or \cite{blocki-lecture-notes},
Lemma 5.3).

\begin{lemma}\label{L:hermit}
Let us fix $0<\lam<\Lam<\infty$. One can find a natural number $N$,
unit vectors $\zeta_1,\dots,\zeta_N\in\HH^n$, and
$0<\lam_*<\Lam_*<\infty$ depending on $\lam,\Lam,n$ only such that
every hyperhermitian $(n\times n)$-matrix $A$ whose spectrum lies in
$[\lam,\Lam]$ can be written
$$\sum_{k=1}^N\beta_k\zeta_k\otimes \zeta_k^* \mbox{ with }
\beta_k\in [\lam_*,\Lam_*].$$ The vectors $\zeta_1,\dots,\zeta_N$
can be chosen to contain any orthonormal basis of $\HH^n$.
\end{lemma}
{\bf Proof.} For convenience of the reader we outline the argument
which is not novel. Let us denote by $K$ the set of quaternionic
hermitian matrices whose spectrum lies in $[\lam,\Lam]$. This is a
compact subset of the interior of the cone of positive definite
hyperhermitian matrices; we denote this open cone by $\cc$. Then
there exists a convex compact polytope $P\subset \cc$ which contains
a neighborhood of $K$. Let $Vert(P)$ denote the set of vertices of
$P$. Using a diagonalization, every matrix $X\in \cc$ can be written
in the form
\begin{eqnarray}\label{E:2-1}
X=\sum_{i=1}^n\alp_i(X) \left(\xi_i(X)\otimes \xi_i^*(X)\right),
\end{eqnarray}
where $\alp_i(X)>0$, and $\xi_i(X)\in\HH^n$ are unit vectors. Let us
fix such a decomposition for any vertex of $P$. Let us define a new
finite subset $S_1\subset \cc$ consisting of rank one non-negative
definite matrices as follows
\begin{eqnarray*}
S_1:=\left\{\left(\sum_j\alp_j(X)\right)\cdot
\left(\xi_i(X)\otimes\xi_i^*(X)\right)|\, X\in Vert(P),
\alp_j(X),\xi_i(X) \mbox{ satisfy } (\ref{E:2-1})\right\}.
\end{eqnarray*}
Let us add to $S_1$ matrices of the form $e_i\otimes e_i^*$,
$i=1,\dots,n$, where $e_1,\dots,e_n$ are an orthonormal basis of
$\HH^n$. Let us denote by $S$ the new set. It is clear that
$P\subset conv(S)$. Hence $conv(S)$ contains a neighborhood of $K$.
Now the lemma follows from the following general fact which is left
to the reader (and where one takes $Q=conv(S)$): Let $K$ be a
compact subset of $\RR^N$ which is contained in the interior of a
compact convex polytope $Q$. Then there exists $\eps>0$ such that
any point $x\in K$ can be written as a convex combination of
vertices of $Q$:
$$x=\sum_v \beta(v)v \mbox{ with } \beta(v)>\eps,$$
where the sum runs over all vertices of $Q$, $\sum_v\beta(v)=1$.
\qed

\hfill

{\bf Proof of Theorem \ref{thm-local}.} The eigenvalues of
$U=(u_{\bar i j})$ belong to $[\lam_1,\Lam_1]$ with
$0<\lam_1<\Lam_1<\infty$ depending on $||\Delta u||_{C^0(\co)},
||f||_{C^0(\co)}$ only. Hence the eigenvalues of $e^fU^{-1}$ are in
$[\lam,\Lam]$ with $0<\lam<\Lam<\infty$ under control. By Lemma
\ref{L:hermit} there exist $N$, unit vectors
$\zeta_1,\dots,\zeta_N\in \HH^n$, and $0<\lam_*<\Lam_*<\infty$ such
that for any $y\in\co$
$$e^{f(y)}U^{-1}(y)=\sum_{i=1}^N\beta_k(y)\zeta_k\otimes \zeta_k^*
\mbox{ with }\beta_k(y)\in [\lam_*,\Lam_*].$$

Observe also that for a unit vector $\zeta\in \HH^n$
$$Tr((\zeta\otimes\zeta^*)(u_{\bar i j}))=Tr(\zeta^*(u_{\bar i
j})\zeta)=\Delta_\zeta u.$$ This and (\ref{E:II})-(\ref{E:II-1})
imply
\begin{eqnarray}\label{E:III}
\sum_{k=1}^N\beta_k(y)(\Delta_{\zeta_k}u(y)-\Delta_{\zeta_k}u(x))\leq
C_1||x-y|| \mbox{ for } x,y\in \co.
\end{eqnarray}

Let $\co''$ be a compact neighborhood of $\co'$ in $\co$. Consider
now a ball $B_{4r}=B(z_0,4r)\subset \co$ with center $z_0\in \co''$.
Let us denote
$$M_{k,r}:=\sup_{B_r}\Delta_{\zeta_k}u,\,
m_{k,r}:=\inf_{B_r}\Delta_{\zeta_k}u,$$
$$\eta(r):=\sum_{k=1}^N(M_{k,r}-m_{k,r}).$$ We will show that for
some $\alp\in (0,1), r_0>0,C>0$ under control
$$\eta(r)\leq Cr^\alp \mbox{ for } 0<r<r_0.$$

Since $\zeta_1,\dots,\zeta_N$ can be chosen to contain an
orthonormal basis of $\HH^n$, this will imply an estimate on
$||\Delta u||_{C^\alp(\co'')}$. Then the Schauder estimates
(\cite{gilbarg-trudinger}, Theorem 4.6) will imply an estimate on
$||u||_{C^{2,\alp}(\co')}$.

The condition $\eta(r)\leq Cr^\alp$ is equivalent to
$$\eta(r)\leq \delta \eta(4r)+r,\, 0<r<r_1$$
where $\delta\in (0,1), r_1$ are under control
(\cite{gilbarg-trudinger}, Lemma 8.23). Summing up
(\ref{harnack-main}) over $\zeta=\zeta_l$ with $l\ne k$ for any
fixed $k$ we get
\begin{eqnarray}\label{E:IV}
r^{-4n}\int_{B_r}\sum_{l\ne k}(M_{l,4r}-\Delta_{\zeta_l}u)\leq
C_3(\eta(4r)-\eta(r)+r).
\end{eqnarray}

But by (\ref{E:III}) we have for any $x\in B_{4r},y\in B_r$
\begin{eqnarray*}
\beta_k(y)(\Delta_{\zeta_k}u(y)-\Delta_{\zeta_k}u(x))\leq
C_1||x-y||+\sum_{l\ne
k}\beta_l(y)(\Delta_{\zeta_l}u(x)-\Delta_{\zeta_l}u(y))\leq\\
C_4r+\Lam_*\sum_{l\ne k}(M_{l,4r}-\Delta_{\zeta_l}u(y)).
\end{eqnarray*}

Optimizing in $x$ we get $$\Delta_{\zeta_k}u(y)-m_{k,4r}\leq
\frac{1}{\lam_*}\left(C_4r+\Lam_*\sum_{l\ne
k}(M_{l,4r}-\Delta_{\zeta_l}u(y))\right).$$ Integrating the last
inequality over $y\in B_r$ and using (\ref{E:IV}) we obtain
\begin{eqnarray}\label{E:V}
r^{-4n}\int_{B_r}(\Delta_{\zeta_k}u(y)-m_{k,4r})\leq
C_5(\eta(4r)-\eta(r)+r).
\end{eqnarray}
Let us estimate the left hand side of (\ref{E:V}) from below. Since
we have normalized the Lebesgue measure on $\HH^n$ so that
$vol(B_1)=1$, we have
\begin{eqnarray*}
r^{-4n}\int_{B_r}(\Delta_{\zeta_k}u(y)-m_{k,4r})\overset{(\ref{harnack-main})}{\geq}
-m_{k,4r}+M_{k,4r}-C(M_{k,4r}-M_{k,r}+r)=\\
C(M_{k,r}-m_{k,r})-(C-1)(M_{k,4r}-m_{k,4r})+C(m_{k,r}-m_{k,4r})-Cr\geq\\
C(M_{k,r}-m_{k,r})-(C-1)(M_{k,4r}-m_{k,4r})-Cr.
\end{eqnarray*}
Substituting this back into (\ref{E:V}) we obtain
\begin{eqnarray*}
C(M_{k,r}-m_{k,r})-(C-1)(M_{k,4r}-m_{k,4r})-Cr\leq
C_5(\eta(4r)-\eta(r)+r).
\end{eqnarray*}
Summing this up over $k$ we get
\begin{eqnarray*}
C\eta(r)-(C-1)\eta(4r)\leq C_6(\eta(4r)-\eta(r)+r).
\end{eqnarray*}
Hence $$\eta(r)\leq \frac{C+C_6-1}{C+C_6}\eta(4r)+r.$$ Theorem
\ref{thm-local} is proved. \qed

\section{Proof of the main theorem.}\label{S:higher-order}
Let us assume that our compact connected hypercomplex manifold $M$
with an HKT-form $\Ome_0$ admits a real (in the quaternionic sense)
and strictly positive (in particular, nowhere vanishing)
$(2n,0)$-form $\Theta\in C^\infty(M,\Lam^{2n,0}_{I,\RR})$ which is
$I$-holomorphic, i.e. $\bar\pt \Theta=0$. Consequently one has
$$\pt\bar\Theta=\pt_J\bar\Theta=0.$$
For any integer $k\geq 1$ and $\beta\in (0,1)$ let us define
\begin{eqnarray*}
\cu^{k,\beta}:=\{\phi\in C^{k,\beta}(M)|\, \Ome_0+\pt\pt_J\phi>0
\mbox{ and } \int_M \phi\cdot
\Ome_0^n\wedge\bar\Ome_0^n=0 \},\\
\cv^{k,\beta}:=\{\chi\in C^{k,\beta}(M)|\, \chi>0 \mbox{ and }
\int_M(\chi-1)\cdot\Ome_0^n\wedge\bar\Theta=0\}.
\end{eqnarray*}
Define
$$\cm(\phi):=\frac{(\Ome_0+\pt\pt_J\phi)^n}{\Ome_0^n}.$$
We claim that $$\cm\colon \cu^{k+2,\beta}\to \cv^{k,\beta}$$ is a
continuous map. The continuity is obvious, the only thing to check
is that for any $\phi\in C^{k+2,\beta}(M)$ one has
$$\int_M (\cm(\phi)-1)\Ome_0^n\wedge \bar\Theta=0.$$
Indeed the left hand side of the last equality is equal to
\begin{eqnarray*}
\int ((\Ome_0+\pt\pt_J\phi)^n-\Ome_0^n)\wedge\bar\Theta=\int
\pt\pt_J\phi\wedge(\sum_{k=0}^{n-1}(\pt\pt_J\phi)^k\wedge
\Ome_0^{n-k-1})\wedge\bar\Theta=\\
\int\phi\wedge \pt\pt_J\left((\sum_{k=0}^{n-1}(\pt\pt_J\phi)^k\wedge
\Ome_0^{n-k-1})\wedge\bar\Theta\right)=0,
\end{eqnarray*}
where in the last equality we have used the Leibnitz rule,
$\pt^2=\pt_J^2=0$, $\pt\pt_J=-\pt_J\pt$, and
$\pt\Ome_0=\pt_J\Ome_0=\pt\bar\Theta=\pt_J\bar\Theta=0$.

\hfill

Next notice that $\cu^{k+2,\beta}$ and $\cv^{k,\beta}$ are open
subsets in Banach spaces with the induced H\"older norms.
\begin{proposition}\label{P:open-image}
Let $(M^{4n},I,J,K)$ be a compact connected hypercomplex manifold
with an HKT-form $\Ome_0$ and which admits a form $\Theta$ as above.
Let $k\geq 1$ be an integer, and let $\beta\in (0,1)$. Then the map
$\cm\colon \cu^{k+2,\beta}\to \cv^{k,\beta}$ is locally a
diffeomorphism of Banach spaces, and in particular its image is an
open subset.
\end{proposition}
{\bf Proof.} By the inverse function theorem for Banach spaces, it
suffices to show that the differential of $\cm$ at any $\phi\in
\cu^{k+2,\beta}$ is an isomorphism of tangent spaces. The tangent
space to $\cu^{k+2,\beta}$ at any point is
$$\tilde C^{k+2,\beta}(M):=\{\phi\in C^{k+2,\beta}(M)|\, \int_M
\phi\cdot\Ome_0^n\wedge\bar\Ome_0^n=0\}.$$ The tangent space to
$\cv^{k,\beta}$ at $\cn(\phi)$ is
$$\tilde C^{k,\beta}(M):=\{\chi\in C^{k,\beta}(M)|\, \int_M \chi
\cdot \Ome_0^n\wedge\bar\Theta =0\}.$$

The differential of $\cm$ at $\phi$ is equal to
\begin{eqnarray*}
D\cm_\phi(\psi)=n\frac{(\Ome_0+\pt\pt_J\phi)^{n-1}\wedge\pt\pt_J\psi}{\Ome^n_0}.
\end{eqnarray*}
Defined by this formula we consider $D\cm_\phi$ as a map
$C^{k+2,\beta}(M)\to C^{k,\beta}(M)$ (without restricting to $\tilde
C$). Then obviously $D\cm_\phi$ is a linear second order
differential elliptic operator. It has no free term (i.e.
$D\cm_\phi(1)=0$), and its coefficients belong to $C^{k,\beta}$.

By the strong maximum principle (see e.g. \cite{gilbarg-trudinger},
Theorem 8.19) and since $M$ is connected, the kernel of $D\cm_\phi$
consists only of constant functions; thus it is one dimensional. The
image of $D\cm_\phi$ is a closed subspace of $C^{k,\beta}(M)$ by
\cite{joyce}, Theorem 1.5.4, which is a version of the Schauder
estimates.

Since the symbol of any second order differential operator of real
valued functions is self-adjoint, the index of $D\cm_\phi$ equals 0
(for operators with coefficients from H\"older spaces see the book
\cite{joyce}, the proof of Theorem C3 in \S 5.6 and Theorem 1.5.4).
Hence $codim\, Im(D\cm_\phi)=1$. But since $Im(D\cm_\phi)\subset
\tilde C^{k,\beta}(M)$, and since $\dim C^{k,\beta}(M)/\tilde
C^{k,\beta}(M)=1$, it follows that $Im(D\cm_\phi)=\tilde
C^{k,\beta}(M)$. \qed

\begin{proposition}\label{P:closed-image}
Let $(M^{4n},I,J,K)$ be a compact manifold with a locally flat
hypercomplex structure which admits a flat hyperK\"ahler metric
compatible with the hypercomplex structure. Let $\Ome_0$ be an
HKT-form (not necessarily corresponding to the hyperK\"ahler
metric). Let $k\geq 2$ be an integer, $\beta\in (0,1)$. Then the map
$$\cm\colon \cu^{k+2,\beta}\to \cv^{k,\beta}$$
is a diffeomorphism of Banach manifolds, in particular it is onto.
\end{proposition}
{\bf Proof.} $\cm$ is one-to-one by the uniqueness of the solution
(in \cite{alesker-verbitsky-10}, Corollary 4.10, the uniqueness was
proven for $C^\infty$-solutions, but exactly the same standard
proof, based on ellipticity and the strong maximum principle, works
under the current assumptions on smoothness).

Now notice that the assumptions of Proposition \ref{P:open-image}
are satisfied. Indeed let $G$ be a locally flat hyperK\"ahler
matric. Let $\Ome$ be the corresponding HKT-form. Then
$\Theta=\Ome^n$ satisfies the assumptions of Proposition
\ref{P:open-image}.

Thus by Proposition \ref{P:open-image} it suffices to show that
$\cm$ is onto. Since $\cv^{k,\beta}$ is obviously connected (it is
even convex), and since the image of $\cm$ is open by Proposition
\ref{P:open-image}, it suffices to show that the image of $\cm$ is a
closed subset of $\cv^{k,\beta}$.

Let we have a sequence of point in the image
$\cm\phi_i\overset{C^{k,\beta}}{\to} e^f\in \cv^{k,\beta}$ where
$\phi_i\in \cu^{k+2,\beta}$. By Theorems \ref{T:second-order},
\ref{thm-main}, and the zero order estimate in
\cite{alesker-verbitsky-10}, Corollary 5.7 (see also
\cite{alesker-shelukhin}, Theorem 2), there exist $\alp\in (0,1)$
and a constant $C$ both depending on $||f||_{C^2}$, $(M,I,J,K)$,
$\Ome_0$, and the locally flat hyperK\"ahler metric, such that
$||\phi_i||_{C^{2,\alp}}<C$ for $i\gg 1$. By the Arzel\`a-Ascoli
theorem choosing a subsequence we may assume that $\phi_i\to \phi$
in $C^2(M)$. Clearly $\phi\in C^{2,\alp}(M)$ and one has
$\cm(\phi)=e^f$. In other words one has
\begin{eqnarray}\label{E:limit-solut}
(\Ome_0+\pt\pt_J\phi)^n=e^f\Ome_0^n.
\end{eqnarray}
Also clearly $\Ome_0+\pt\pt_J\phi\geq 0$. But because of
(\ref{E:limit-solut}) the inequality is strict, i.e.
$\Ome_0+\pt\pt_J\phi>0$, and the equation (\ref{E:limit-solut}) is
elliptic with $C^\infty$ coefficients on the left hand side and with
$C^{k,\beta}$ on the right hand side. Hence by Lemma 17.16 from
\cite{gilbarg-trudinger} $\phi\in C^{k+2,\beta}(M)$. Thus $\phi\in
\cv^{k+2,\beta}$ and $e^f=\cm(\phi)\in Im(\cm)$. \qed

\hfill

Finally let us state the main result of the paper which is an
immediate consequence of Proposition \ref{P:closed-image}.
\begin{theorem}\label{T:MA-solution}
Let $(M^{4n},I,J,K)$ be a compact locally flat hypercomplex manifold
which admits a flat hyperK\"ahler form (of class $C^\infty$)
compatible with the hypercomplex structure. Let $\Ome_0$ be an
HKT-form on $M$ of class $C^\infty$ (which does not necessarily
correspond to the above mentioned hyperK\"ahler metric). Let $k\geq
2$ be an integer, $\beta\in (0,1)$. Let $f\in C^{k,\beta}(M)$.

Then there is a unique constant $A>0$ such that quaternionic
Monge-Amp\`ere equation
\begin{eqnarray}\label{E:MA-in-main-thm}
(\Ome_0+\pt\pt_J\phi)^n= A e^f\Ome_0^n
\end{eqnarray}
has a unique, up to a constant, $C^2$ smooth solution $\phi$ which
necessarily belongs to $C^{k+2,\beta}(M)$. If $f$ is $C^\infty$
smooth, then any solution $\phi$ is also $C^\infty$ smooth.
\end{theorem}

\begin{remark}
The constant $A$ in the theorem is defined by the following
condition. Let $\Ome$ be the HKT-form corresponding to the locally
flat hyperK\"ahler metric. Then $A$ is found from the equality
$$\int_M
Ae^f\cdot \Ome_0^n\wedge\bar\Ome^n=\int_M\Ome_0^n\wedge\bar\Ome^n.$$
\end{remark}

\end{document}